# A Method of Advanced Chaos Optimal Search Algorithm

UnSun Pak, YongNam Kim, GyongIl Ryang

Faculty of Electronics & Automation, **Kim Il Sung** University, Pyongyang,

Democratic People's Republic of Korea

**Abstract:** In this paper, the advanced parallel chaos optimal search algorithm is proposed and the effectiveness of the proposed algorithm is verified through the experiment for find out the minimum of several benchmark functions.

**Keywords**: chaos optimal search, chaos search, optimal search, search algorithm

## 1. Introduction

Chaos theory has been applied to many branches including encryption, communication, graphic data compression, control and optimal search, etc.

And especially, the researches for use chaos to optimal search have been intensified.

In previous work [1], the chaos optimal search algorithm that does not use the gradient of function was proposed, but the drawback is it's long searching time.

To overcome this drawback, variable scaling chaos optimal search algorithm[2,3], chaos optimal search algorithm using variable reduction rate[4,5], hybrid chaos optimal search algorithm with simulated annealing method and GA[6,7] was proposed. But the drawback of the variable scaling chaos optimal search algorithm and the chaos optimal search algorithm using variable reduction is that it could fall into a local optimal point, and the drawback of the hybrid chaos optimal search algorithm is that it is complicated to accomplish.

In this paper, to overcome the drawback of the previous chaos optimal search algorithms, we propose the advanced parallel chaos optimal search algorithm and verify the effectiveness of the proposed algorithm through the optimal search experiment for several benchmark functions.

## 2. Advanced parallel chaos optimal search algorithm(APCOSA)

Function optimization problem is formulated as follows.

$$\min_{x} f(x) \qquad (1)$$

Here, $x = [x_1, x_2, \cdots, x_n]^T$ is decision variable vector and $a_i \le x_i \le b_i, i = 1,2,\cdots,n$.

The APCOSA consist of several stages including initialization, rough search, elimination of the neighboring points, parallel search and fine search. And the search parameters are the count of inner iteration $N$, the count of outer iteration $M$, the count of the initial candidate points $p$.

**[Stage 1 : Initialization]**

The candidate point set $X_{\min} = \{x_{\min}^1, x_{\min}^2, \cdots, x_{\min}^p\}$ is composed of $p$ points which randomly generated in search space and the function value set $F_{\min} = \{f_{\min}^1, f_{\min}^2, \cdots, f_{\min}^p\}$ which evaluated on the candidate point set is obtained.

**[Stage 2 : Rough Search]**

The initial point $X^0 = [X_1^0, X_2^0, \cdots, X_n^0]^T$ is randomly generated in range $[0, 1]^n$ and from this point; the point sequence $\{X^k, k = 0,1,2,\cdots\}$ is generated by the Logistic mapping.

$$X_i^{k+1} = AX_i^k(1-X_i^k), \quad i=1,2,\cdots,n, \quad k=0,1,2,\cdots \tag{2}$$

Then variable $X^k$ is ergodic in the range $[0,1]^n$ [2].

$X^k = [X_1^k, X_2^k, \cdots, X_n^k]^T$ is scaled to the defining range and the decision variable vector $x^k = [x_1^k, x_2^k, \cdots, x_n^k]^T$ is obtained. That is

$$x_i^k = a_i + (b_i - a_i) \cdot X_i^k \tag{3}$$

In the candidate point set $X_{\min} = \{x_{\min}^1, x_{\min}^2, \cdots, x_{\min}^p\}$, the point $x_{\min}^l, 1 \le l \le p$ which nearest to the point $x^k = [x_1^k, x_2^k, \cdots, x_n^k]^T$ is found. If $f(x^k) < f_{\min}^l$ then $x_{\min}^l$ and $f_{\min}^l$ is substituted for $x^k$ and $f(x^k)$.

This process is iterated for $k = 1, 2, \cdots, N$.

**[Stage 3 : Elimination of the Neighboring Points]**

After the rough search, neighboring points could be exist in the candidate point set $X_{\min} = \{x_{\min}^1, x_{\min}^2, \cdots, x_{\min}^p\}$. Without eliminate these points, search is duplicated in the next stage, the parallel search stage, and the searching effectiveness is fallen down. So in this stage, neighboring points except only one point are eliminated from the candidate point set.

First, neighboring distance $d_{near}$ is defined as follows.

$$d_{near} = \frac{\sqrt{\sum_{i=1}^{n}(b_i - a_i)^2}}{2p} \tag{4}$$

Next, two points that the distance between them is smaller than $d_{near}$ is found out and one of them is eliminated. This operation is iterated until any distance between two existing points is no smaller than $d_{near}$. Finally new sets $X'_{\min} = \{x_{\min}^1, x_{\min}^2, \cdots, x_{\min}^{p'}\}$ and $F'_{\min} = \{f_{\min}^1, f_{\min}^2, \cdots, f_{\min}^{p'}\}$ is obtained.

**[Stage 4 : Parallel Search]**

Every point of the candidate point set $X'_{\min} = \{x_{\min}^1, x_{\min}^2, \cdots, x_{\min}^{p'}\}$ has become the central point and the parallel search is done.

Local chaos search at the every candidate point $x_{\min}^l$ consist of $M-2$ substage and in the $m(\le M-2)^{th}$ substage the searching range is reduced as follows.

$$\begin{cases} a_i^m = \max\left\{a_i, \quad x_{\min,i}^l - (b_i - a_i)\dfrac{1}{2p' \cdot m!}\right\}, \\ b_i^m = \min\left\{b_i, \quad x_{\min,i}^l + (b_i - a_i)\dfrac{1}{2p' \cdot m!}\right\}, \quad m=1,2,\cdots,M-2 \end{cases} \tag{5}$$

Chaos search in every substage is done as follows.

The initial point $X^0 = [X_1^0, X_2^0, \cdots, X_n^0]^T$ is randomly generated in range $[0,1]^n$ and from this point; the point sequence $\{X^k, k=0,1,2,\cdots,N'\}$ is generated by the Logistic mapping.(2). Here the count of iteration is $N' = [N/p']$. $X^k = [X_1^k, X_2^k, \cdots, X_n^k]^T$ is scaled to the reduced searching range and the decision variable vector $x^k = [x_1^k, x_2^k, \cdots, x_n^k]^T$ is obtained. If $f(x^k) < f_{\min}^l$ then $x_{\min}^l$ and $f_{\min}^l$ is substituted for $x^k$ and $f(x^k)$.

This process is iterated for $k = 1, 2, \cdots, N'$.

The candidate point set is continuously updated during the parallel search. After the parallel search, the best point in the candidate point set is denoted $x_{min}$.

**[Stage 5 : Fine Search]**

Here, for the obtained point $x_{min}$, chaos search is done with the searching range

$$a'_i = \max\left\{a_i, x_{min} - (b_i - a_i)\frac{1}{2p' \cdot (M-1)!}\right\},$$
$$b'_i = \min\left\{b_i, x_{min} + (b_i - a_i)\frac{1}{2p' \cdot (M-1)!}\right\} \quad (6)$$

and the count of iteration $N$.

Finally, optimal solution $(x^*, f^*)$ is obtained.

For APCOSA, the evaluation count of function is $M \cdot N$ and the tolerance of obtained solution is $(b_i - a_i)/M!$ for the $i^{th}$ variable.

## 3. Experiment

The effectiveness of the proposed algorithm is verified through the experiment for find out the minimum of follow 5 benchmark functions.

$$F_1(x_1, x_2) = 100(x_1^2 - x_2)^2 + (1 - x_1)^2, \quad -2.048 \leq x_1, x_2 \leq 2.048$$
$$F_2(x_1, x_2) = (4 - 2.1x_1^2 + \frac{x_1^4}{3})x_1^2 + x_1 x_2 + (-4 + 4x_2^2)x_2^2, \quad -100 < x_1, x_2 < 100$$
$$F_3(x_1, x_2) = -0.5 + \frac{\sin^2\left(\sqrt{x_1^2 + x_2^2}\right) - 0.5}{\left(1 + 0.001(x_1^2 + x_2^2)\right)^2}, \quad -100 < x_1, x_2 < 100$$
$$F_4(x_1, x_2) = (x_1^2 + x_2^2)^{0.25}\left[\sin^2(50(x_1^2 + x_2^2)^{0.1} + 1)\right], \quad -100 < x_1, x_2 < 100 \quad (7)$$
$$F_5(x_1, x_2) = \left[1 + (x_1 + x_2 + 1)^2(19 - 14x_1 + 3x_1^2 - 14x_2 + 6x_1 x_2 + 3x_2^2)\right]$$
$$\times \left[30 + (2x_1 - 3x_2)^2(18 - 32x_1 + 12x_1^2 + 48x_2 - 36x_1 x_2 + 27x_2^2)\right],$$
$$-2 \leq x_1, x_2 \leq 2$$

The graph of these benchmark functions is shown fig.1.

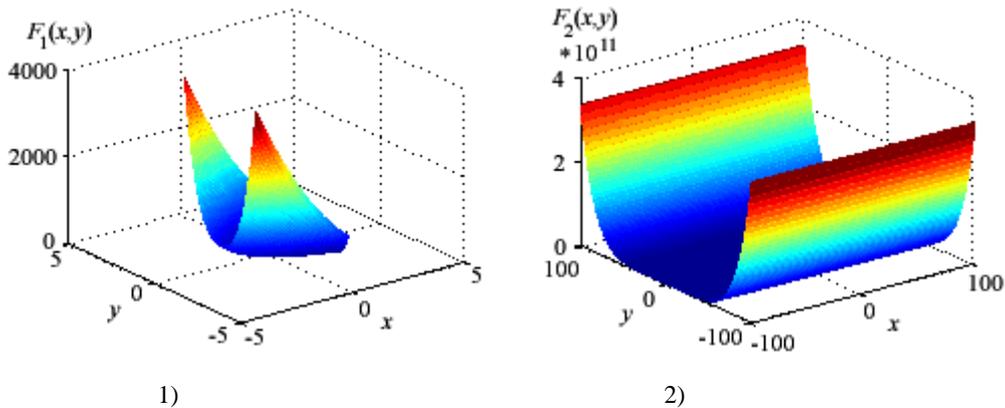

1)　　　　　　　　　　　　　　　　2)

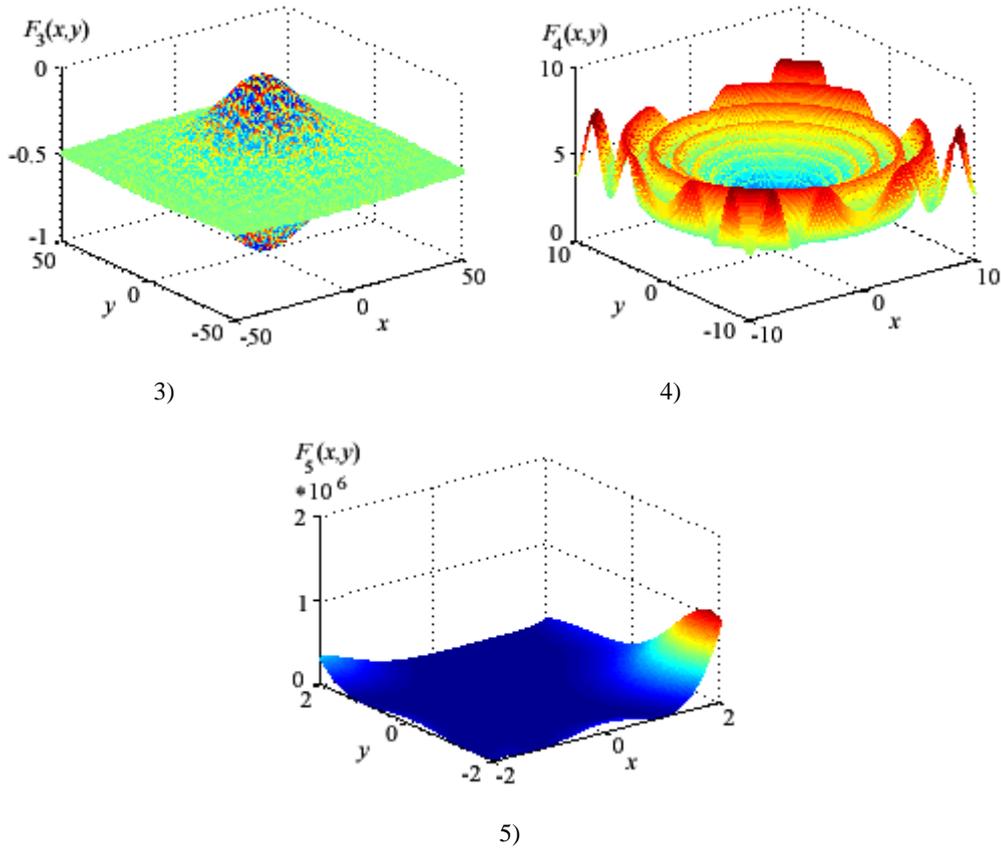

Fig.1. the graph of benchmark functions

1) $F_1(x, y)$, 2) $F_2(x, y)$, 3) $F_3(x, y)$, 4) $F_4(x, y)$, 5) $F_5(x, y)$

In experiment, the effectiveness of the proposed algorithm compare with the chaos optimal search algorithm[1], chaos optimal search algorithm using variable reduction rate[4]. The result is shown table 1.

Table1. The result of experiment
(1-Chaos optimal search algorithm, 2-Chaos optimal search algorithm using variable reduction rate, 3-APCOSA)

| № | function | minimum $(x_1^*, x_2^*)$ | $F^*$ | Searching time(s) 1 | 2 | 3 | Success frequency 1 | 2 | 3 |
|---|---|---|---|---|---|---|---|---|---|
| 1 | $F_1$ | (1, 1) | 0 | 0.9257 | 0.3600 | 0.3560 | 89/100 | 99/100 | 100/100 |
| 2 | $F_2$ | (0.0898, -0.7127) | -1.0316 | 1.356 | 0.6210 | 0.6020 | 85/100 | 95/100 | 99/100 |
| 3 | $F_3$ | (0, 0) | 1 | 3.793 | 1.682 | 1.589 | 70/100 | 82/100 | 98/100 |
| 4 | $F_4$ | (0, 0) | 0 | 4.312 | 3.335 | 3.305 | 85/100 | 95/100 | 99/100 |
| 5 | $F_5$ | (0, -1) | 3 | 3.973 | 1.642 | 1.565 | 89/100 | 99/100 | 100/100 |

As shown in the table, APCOSA has shorter searching time and higher success frequency compare with the previous algorithms.

## 4. Conclusion

First, the advanced parallel chaos optimal search algorithm had been proposed.

Second, the effectiveness of proposed algorithm is verified through the comparison with the previous chaos optimal search algorithms.